\documentclass[12pt]{amsart}
\usepackage{graphicx}
\usepackage{amssymb}
\usepackage{amsfonts}
\usepackage{amsmath}
\usepackage{array}
\usepackage{rotating}

\newtheorem{example}{Example}[section]

\newtheorem{theorem}[example]{Theorem}

\newtheorem{corollary}[example]{Corollary}

\newtheorem{proposition}[example]{Proposition}

\def\Proof{\noindent \it Proof -- \rm}
\def\qed{\hspace{3.5mm} \hfill \vbox{\hrule height 3pt depth 2 pt width 2mm}
\bigskip}

\def\FQSym{{\bf FQSym}}

\def\Sym{{\bf Sym}}
\def\NCSF{{\bf Sym}}

\def\<{\langle}
\def\>{\rangle}

\def\Z{\operatorname{\mathbb Z}}

\def\F{{\bf F}}

\def\G{{\bf G}}

\def\SG{{\mathfrak S}}

\def\maj{{\rm maj}}

\def\std{{\rm std}}

\def\DSym{{\bf DSym}}
\def\DQSym{{\bf DQSym}}
\def\id{{\rm id}}

\def\kRC{{\rm RC}}
\def\kDC{{\rm DC}}

\def\NkRC{{\rm N}}
\def\kR{{\sf R}}
\def\kS{{\sf S}}
\def\kF{{\sf F}}

\def\kE{{\sf E}}

\def\SF{{\bf S}}
\def\EF{{\bf E}}

\title[Generalized descent patterns and Hopf algebras]%
{Generalized descent patterns in permutations\\
and associated  Hopf algebras}

\author[J.-C.~Novelli,  C. Reutenauer, J.-Y.~Thibon]%
{Jean-Christophe Novelli, Christophe Reutenauer, Jean-Yves Thibon}

\address[Novelli, Thibon] {Institut Gaspard Monge, Universit\'e Paris-Est
Marne-la-Vall\'ee, \\ 5 Boulevard Descartes \\ Champs-sur-Marne \\
77454 Marne-la-Vall\'ee cedex 2 \\ France}
\address[Reutenauer]{LaCIM,
Universit\'e du Qu\'ebec \`a Montr\'eal\\
CP 8888, Succ. Centre-ville\\
Montr\'eal (Qu\'ebec) H3C 3P8}
\email[Jean-Christophe Novelli]{novelli@univ-mlv.fr (corresponding author)}
\email[Christophe Reutenauer]{Reutenauer.Christophe@uqam.ca}
\email[Jean-Yves Thibon]{jyt@univ-mlv.fr} 
\date{\today}

\begin{document}

\begin{abstract}
Descents in permutations or words are defined from the relative position of
two consecutive letters. We investigate a statistic involving patterns
of $k$ consecutive letters, 
and show that it leads to Hopf algebras generalizing noncommutative symmetric
functions and quasi-symmetric functions.
\end{abstract}

\maketitle

\section{Introduction and Background}

Recall that the \emph{standardized word} $\std(w)$ of a word $w\in A^*$
over an ordered alphabet $A$
is the permutation obtained by iteratively scanning $w$ from left to right,
and labelling $1,2,\ldots$ the occurrences of its smallest letter, then
numbering the occurrences of the next one, and so on.
For example, $\std(bbacab)=341625$.

A permutation $\sigma\in\SG_n$ is said to have a descent at $i$ if
$\sigma_i>\sigma_{i+1}$. One can alternatively say that the
standardization of the two-letter word $\sigma_i\sigma_{i+1}$ is
$21$, and the descent set of $\sigma$ can be encoded by the
descent pattern $(\std(\sigma_i\sigma_{i+1})_{i=1\dots n-1})$,
produced by scanning $\sigma$ with a sliding window of width two.

An obvious generalization of this notion would be to use a window of arbitrary
width $k$. For example, the permutation $\sigma=(85736124)$ would have as
3-descent pattern the sequence
\begin{equation}
\label{fenetre}
p = (312, 231, 312, 231, 312, 123) \,.
\end{equation} 

This idea immediately raises a couple of questions. It is known that the sums
of permutations of $\SG_n$ having the same descent pattern span a subalgebra
$\Sigma_n$ of the group algebra (Solomon's descent algebra \cite{So}) and that the
direct sum of all the $\Sigma_n$ has a natural Hopf algebra structure
(Noncommutative symmetric functions), inherited from that of the Hopf algebra
of permutations~\cite{MR1,NCSF6}.
Are there analogs of these facts for the generalized descent classes?

We will show that, although the $k$-descent classes do not span a subalgebra
of the group algebra, the Hopf algebra construction still works, and leads
to some interesting combinatorics.

\section{Generalized descent patterns and codes}

For ordinary descents in $\SG_n$, a classical encoding of the patterns is by
compositions of $n$.
Recall that if the descents of $\sigma$ form the set
$D=\{d_1,\ldots,d_{r-1}\}$, we encode it by the composition
$I=(i_1,\ldots,i_r)$ of $n$ such that $d_j=i_1+i_2+\cdots i_j$. Then, in the
algebra of noncommutative symmetric functions, the concatenation of
compositions corresponds to the (outer) multiplication in various bases,
appropriately called \emph{multiplicative bases}.

\subsection{The $k$-descent code}

Another encoding having the same property would be to represent the
descent set $D$ by a binary word $b=b_1b_2\cdots b_n$, with $b_i=1$
if $i\in D$ and $b_i=0$ otherwise (note that $b_n$ is always 0).
It is this kind of encoding which is most easily generalized to
$k$-descent patterns.

\noindent
We regard $\sigma\in\SG_n$ as a permutation of $\Z$ whose support is contained
in $[1,n]$ and associate with it the sequence $(d_i\in[1,k])_{i\in\Z}$, where
$p_i$ is the relative position of $\sigma(i)$ w.r.t. its $k-1$ predecessors
$\sigma(i-1),\ldots,\sigma(i-k+1)$, that is, $d_i=j$ if $\sigma(i)$ is the
$j$th element of the sequence
$\sigma(i),\sigma(i-1),\ldots,\sigma(i-k)$ sorted in increasing order.
Hence $d_i=k$ for $i\le 1$ and for $i>n$. 

We can therefore identify the $k$-descent pattern with the word of length $n$
\begin{equation}
\kDC_k(\sigma) = d_1d_2\cdots d_n \in [1,k]^n\,,
\end{equation}
which we will call the $k$-descent code.
Indeed, starting from the $k$-descent pattern $p=(p_1,\dots,p_r)$ of $\sigma$,
one recovers the $k$-descent code by first computing the $k$-descent code of
$p_1$ and appending to it the final letters of all words $p_2$ up to
$p_r$.
Conversely, $d_1\dots d_k$ gives back $p_1$ since, in that case the sequence
$d_1\dots d_k$ represents the complement to $k$ of the code of $p_1$ 
(obtained by changing $i$ into $k-i$),
and then gives back $p_k$ from both $p_{k-1}$ and $d_k$ since $p_k$ is the
word obtained by removing the first letter of $p_{k-1}$, then adding $1$ to
all values greater than or equal to $d_k$, and concatenating $d_k$ to the 
resulting word.

For example, one can check on
$\kDC_3(85736124)= 32212123$ that the
above algorithm gives back the sequence (\ref{fenetre}).

Here are two more examples of the $k$-descent code.
\begin{equation}
\kDC_3(426135)= 323123 \qquad\text{and}\qquad
\kDC_4(426135)= 434133.
\end{equation}

\subsection{The $k$-recoil code and equivalence classes}

The $k$-descent code of the inverse permutation will be called
the $k$-\emph{recoil code}:
$\kRC_k(\sigma):=\kDC_k(\sigma^{-1})$.
It can be computed without inverting permutations in the following way: first
see $\sigma$ as a permutation of $\Z$ as before. Then, for all $i$, restrict
$\sigma$ to the \emph{values} between $i$ and $i-k+1$. And then associate with
$i$ the position of $i$ in this new word.
For example, $\kRC_3(425163)=323123$ which is coherent with the previous
example since $425163^{-1}=426135$.

Since we only need to compare letter $i$ with the $k$ preceding letters in the
lexicographic order, we can rephrase this construction with the help of the
standardization process.
Indeed, two permutations $\sigma$ and $\tau$ have same $k$-recoil code iff
\begin{equation}
\forall i\leq n-k+1,\
\std(\sigma_{|_{[i,i+k-1]}}) =
\std(\tau_{|_{[i,i+k-1]}}),
\end{equation}
where $\sigma_{|_{[a,b]}}$ means the restriction of $\sigma$ to its
\emph{values} in the interval $[a,b]$.
We shall then write $\sigma\equiv_k\tau$. 

This can be extended to words over an ordered alphabet: we set $u\equiv_kv$
iff $u$ is a rearrangement of $v$ and $\std(u)\equiv_k \std(v)$.
For $k=2$, this is the hypoplactic congruence (see~\cite{NCSF4,Nov}).

For example, with $k=3$, each equivalence class in $\SG_n$ with $n \leq 3$
has one element, and there are $18$ classes in $\SG_4$, among which $6$
non-singleton classes:
\begin{equation}
\label{ex-nonsingl}
[1423,4123],\ [1432,4132],\ [2143,2413],\
[2314,2341],\ [3142,3412],\ [3214,3241].
\end{equation}
For example, the $3$-recoil code of both $2314$ and $2341$ is $3223$.

\bigskip
Note that the first letter of the $k$-descent (or recoil) code is always $k$
and the next one is either $k$ or $k-1$.
More generally, a $k$-recoil code is a word $I=(i_1,\dots,i_r)$ satisfying
$i_\ell \in [\max(k-\ell+1,1),k]$ for all $\ell$.
Conversely, given a word $I$ satisfying these conditions, one easily builds a
permutation with $I$ as $k$-recoil code. By induction, there exists a
permutation $\sigma$ with $k$-recoil code $I'=(i_1,\dots,i_{r-1})$. Now, place
$r$ anywhere between the $i_r$-th and the $1+i_r$-th element of $\sigma$ in
the interval $[r-k+1,r-1]$. This permutation has $I$ as $k$-recoil code.
Note that in particular this allows one to build easily the smallest (resp.
the largest) elements for the lexicographic order of each equivalence class:
put at each step letter $r$ at the rightmost (resp. leftmost) possible
spot.

\begin{proposition}
The number $\NkRC(k,n)$ of $k$-descent (or recoil) classes of $\SG_n$ is
\begin{equation}
\left\{
\begin{array}{lr}
n!  & \text{if $n\leq k$,} \\
k!\, k^{n-k}  & \text{if $n\geq k$.} \\
\end{array}
\right.
\end{equation}
More precisely, the $k$-descent (or recoil) codes $I=(i_1,\dots,i_r)$ satisfy
\begin{equation}
\left\{
\begin{array}{lr}
k-\ell+1 \leq i_\ell \leq k & \text{if $\ell\leq k$,} \\
1        \leq i_\ell \leq k & \text{if $\ell\geq k$.}
\end{array}
\right.
\Longleftrightarrow
i_\ell \in [\max(k-\ell+1,1),k].
\end{equation}
\qed
\end{proposition}

For example, the $3$-recoil codes of all permutations of $\SG_3$ and $\SG_4$
are, taking the permutations in lexicographic order:
\begin{equation}
\label{ex3-codeperms}
333, 332, 323, 322, 331, 321.
\end{equation}
\begin{equation}
\label{ex4-codeperms}
\begin{split}
3333, 3332, 3323, 3322,
3331, 3321, 3233, 3232,
3223, 3223, 3232, 3222,\\
3313, 3312, 3213, 3213,
3312, 3212, 3331, 3321,
3231, 3221, 3311, 3211.
\end{split}
\end{equation}
In particular, one can check that the codes
$3331$, $3321$, $3232$, $3223$, $3312$, and $3213$ occur twice, 
and correspond to the six non-singleton $3$-classes in $\SG_4$ (see
Equation~(\ref{ex-nonsingl})).

\subsection{Classes of permutations having the same $k$-recoil code}

The following proposition generalizes the fact that two permutations have same 
recoil code ($k=2$) iff one can go from one to the other by succesively 
exchanging non-consecutive adjacent values.

\begin{proposition}
\label{prop-ech}
Two permutations $\sigma$ and $\mu$ have same $k$-recoil code iff
one can go from $\sigma$ to $\tau$ by successively exchanging adjacent values
whose difference is at least~$k$.
\end{proposition}

\Proof
Let us write $\sigma\equiv'_k\tau$ if one can go from $\sigma$ to $\tau$ by
exchanging adjacent values whose difference is at least $k$.
Then it is obvious that
\begin{equation}
\sigma\equiv'_k\tau \Rightarrow \sigma\equiv_k\tau.
\end{equation}
So each $\equiv'_k$ class is contained in an $\equiv_k$ class.
In particular, the number of $\equiv'_k$ classes is at least equal to the
number of $\equiv_k$ classes, and if those numbers are equal, then so are the
equivalence classes. 

Now, each $\equiv'_k$ class has a minimal element for the lexicographic
order and this element has no consecutive letters $i$ and $j$ so that
$i-j\geq k$. Let us denote by $W(k,n)$ the set of words such that there are no
consecutive letters $i$ and $j$ so that $i-j\geq k$.

There are at most as many $\equiv'_k$ classes as elements in $W(k,n)$.
Observe now that if one removes $n$ from any word of $W(k,n)$, one obtains
a word in $W(k,n-1)$. Moreover, given a word of $W(k,n-1)$, in order to get a
word in $W(k,n)$, one can put $n$ at $n$ different spots if $n\leq k$ or at
$k$ different spots (after $n-k+1$, $\dots$, $n-1$ or at the end)
if $n\geq k$.
So $|W(k,n)|$ is equal to $\NkRC(k,n)$, the number of $\equiv_k$ classes thanks to
the formula relating $\NkRC(k,n)$ and $\NkRC(k,n-1)$.
\qed

This argument unravels a simple characterization of the minimal element of an
equivalence class: 

\begin{corollary}
\label{cor-minmax}
The minimal elements of the classes of $\equiv_k$ 
are the permutations $\sigma$ such that
no two adjacent letters satisfy $\sigma_i-\sigma_{i+1}\geq k$.

By symmetry, the maximal elements of $\equiv_k$ are the permutations $\sigma$
such that  no two adjacent letters
satisfy $\sigma_{i+1}-\sigma_i\geq k$.

Moreover, the set of maximal elements is obtained from the set of minimal
elements of $\SG_n$ by the transformation
$\tau':=(n+1-\tau_1,\dots,n+1-\tau_n)$.
\qed
\end{corollary}

The proposition implies that the $\equiv_k$ classes split the right
permutohedron into connected components. We can be more precise:

\begin{proposition}
The set of permutations having a given $k$-recoil (respectively descent) code
is  an interval of the right (resp. left) weak order.
\end{proposition}

\Proof
Let $C$ be a $k$-recoil class and let $\alpha_C$ (respectively $\omega_C$) be
its minimal (resp. maximal) element.
If $\sigma$ has same $k$-RC as $\alpha_C$ (and $\omega_C$), then, thanks to
Proposition~\ref{prop-ech}, if $\sigma$ is not minimal, $\sigma$ has a pair
of adjacent letters $\sigma_i$ and $\sigma_{i+1}$ such that
$\sigma_i-\sigma_{i+1}\geq k$. Then one can exchange those two letters and
iterate the process until one reaches the minimal element, so that
$\sigma\geq \alpha_C$ for the right permutohedron. The same argument proves
that $\sigma\leq \omega_C$. 

Conversely, all permutations of the interval $[\alpha_C,\omega_C]$ of the
right permutohedron have same $k$-RC: their inversion sets are contained in
the inversion set of $\omega_C$ and contain the inversion set of $\alpha_C$,
so that letters $i$, $i+1$, $\dots$, $i+k-1$ have same relative positions.
\qed

As in the case of ordinary descents, the order ideals defined by maximal
elements are unions of classes. This property is essential for defining
multiplicative bases in the associated Hopf algebras.

\begin{proposition}
\label{prop-unionintervs}
Let $\omega_C$ be the maximal element of an $\equiv_k$ class.
Then the interval $[\id,\omega_C]$ of the right permutohedron is an union of
$\equiv_k$ classes, where $\id$ is the identity permutation.

Moreover, the interval $[\alpha_C,\omega]$ of the right permutohedron is also
an union of $\equiv_k$ classes, where $\alpha_C$ is the minimal element of an
$\equiv_k$ class and $\omega$ is the maximal permutation.
\end{proposition}

\Proof
By Corollary~\ref{cor-minmax}, the second statement is equivalent to the first
one.

Thus, we must prove that if $x\leq \omega_C$ then the maximal element
$\omega_{C'}$ of its $\equiv_k$ class satisfies $\omega_{C'}\leq \omega_C$.
If $x$ is maximal, we are done. Otherwise, thanks to the characterization of
the maximal elements, we have
\begin{equation}
x = \dots i\,j\dots,
\end{equation}
where $j-i\geq k$.
We shall prove that the permutation $x'$ obtained from $x$ by exchanging $i$
and $j$ also satisfies $x'\leq \omega_C$. Then, since all classes are
intervals, we see by induction on the distance from $x$ to $\omega_{C'}$ that
$\omega_{C'}$ also satisfies $\omega_{C'}\leq \omega_C$.
Consider the subset of the elements of the permutohedron greater than
$x$ such that $i$ and $j$ are not exchanged, that is the set of elements
greater than $x$ and not greater than~$x'$. This set does not contain any
maximal element: the values between $i$ and $j$ in such permutations can only
be either smaller than $i$ or greater than $j$ so that there are always two
consecutive values with difference at least $k$.
\qed

\subsection{$k$-Eulerian numbers and polynomials}

The classical Eulerian polynomials count permutations by their number of
descents, or equivalently, by the number of 2s in their 
2-descent code. 

In this form, the definition can be easily generalized.
We define the $k$-Eulerian polynomial
$E_{n,k}$ as the sum over $\SG_n$ of the product of $t_i$ where $i$ runs
through all entries except the first one of their $k$-recoil code.

For example, with $k=3$, we obtain from~(\ref{ex3-codeperms})
and~(\ref{ex4-codeperms})
\begin{equation}
E_{1,3} = 1;\qquad
E_{2,3} = t_2+t_3;\qquad
E_{3,3} = t_1t_2 + t_1t_3 + 2t_2t_3 + t_2^2 + t_3^2;
\end{equation}
\begin{equation}
E_{4,3} = t_1^2t_2 + t_1^2t_3 + 2\,t_1t_2^2 + 7\,t_1t_2t_3 +
3\,t_1t_3^2 + t_2^3 + 5\,t_2^2t_3 + 3\,t_2t_3^2 + t_3^3.
\end{equation}

\subsection{$k$-Major index}

The classical major index is the sum of the
positions of the descents. We can replace the monomial
$q^{\maj(\sigma)}$ by the product
$\prod q_{C_i}^{i-1}$ where
$C=\kRC(\sigma)$.

For example, with $k=3$, the $k$-major index
polynomials of the first symmetric groups are
\begin{equation}
M_{1,3} = 1;\qquad
M_{2,3} = q_2+q_3;\qquad
M_{3,3} = q_1^2q_2 + q_1^2q_3 + q_2^2q_3 + q_2^3 + q_2q_3^2 + q_3^3.
\end{equation}

\section{Associated combinatorial Hopf algebras}

The equivalence relation $\equiv_k$ can be used to define subalgebras
and quotients of the Hopf algebra of permutations, generalizing respectively
noncommutative symmetric functions and quasi-symmetric functions (cf.
\cite{NCSF1}). 

Recall that the Hopf algebra of permutations, introduced in
\cite{MR1}, can be realized as the algebra $\FQSym$ (Free
Quasi-Symmetric functions, cf. \cite{NCSF6}), spanned by
the polynomials
\begin{equation}
\G_\sigma(A)=\sum_{\std(w)=\sigma}w\,.
\end{equation}

\subsection{Subalgebras}

Imitating the case $k=2$, we define generalized ribbons by
\begin{equation}
\kR_C=\sum_{\kDC_k(\sigma)=C}\G_\sigma
\end{equation}
for a $k$-descent code $C$.

This basis generalizes the classical (strict) descent classes. We can also
generalize the large descent classes (permutations whose descent set
is contained in a prescribed one), corresponding to the multiplicative basis
$S^I$ of $\Sym$. We set
\begin{equation}
\kS^C := \SF^{\omega_C} = \sum_{\sigma\le\omega_C}\G_\sigma
\quad{\rm and} \quad
\kE^C := \EF^{\alpha_C} = \sum_{\sigma\ge\alpha_C}\G_\sigma,
\end{equation}
where $\leq$ and $\geq$ correspond to the left weak order.

The products of the $\SF^\sigma$ and the $\EF^\sigma$ in $\FQSym$ are
well-known \cite{NCSF7}:
\begin{equation}
\SF^\sigma \SF^\tau = \SF^{\sigma[|\tau|]\cdot \tau},
\qquad
\EF^\sigma \EF^\tau = \EF^{\sigma\cdot \tau[|\sigma|]},
\end{equation}
where $\mu[i]=(\mu_1+i,\dots,\mu_n+i)$.

We can then state:

\begin{theorem}
The $\kS^C$ and the $\kE^C$ are multiplicative bases of a subalgebra
$\DSym^{(k)}$ of $\FQSym$. The $\kR_C$ are also a basis of $\DSym^{(k)}$.

Moreover, $\DSym^{(k)}$ is free as an algebra over the $\SF^\sigma$
(respectively $\EF^\sigma$) indexed by permutations that are both
mirror images of connected permutations and maximal (resp. minimal) elements
of $\equiv_k$ classes.
%
%
\end{theorem}

\Proof
Since the shifted concatenation $\sigma\cdot \tau[|\sigma|]$ of two minimal
elements of an $\equiv_k$ class is a minimal element, the product of two
$\kE^C$ is internal, so $\DSym^{(k)}$ is a subalgebra of $\FQSym$ and the
$\kE^C$ are a multiplicative basis. By inclusion-exclusion, we get that the
$\kR_C$ (and the $\kS^C$ too) are bases of $\DSym^{(k)}$. The fact that the
$\SF$ are multiplicative imply that the $\kS$ are also multiplicative.


Moreover, since the $\EF^\sigma$ indexed by permutations that are
mirror images of connected permutations generated $\FQSym$, the free
subalgebra $\DSym^{(k)}$ is free over a set of generators indexed by
permutations that are both mirror images of connected permutations and minimal
elements of $\equiv_k$ classes.
The same holds for the $\SF^\sigma$.
%
\qed

For example,
\begin{equation}
\kR_{321}\kR_{3321} =
\kR_{3211221} + \kR_{3211321} + \kR_{3212321} + \kR_{3213321}\,.
\end{equation}

\medskip

The Hilbert series of $\DSym^{(k)}$ is
\begin{equation}
H_k(t)=\sum_{j=0}^{k-1}j!t^j+\frac{k!t^k}{1-kt}
\end{equation}
and the generating series $G_k(t)$ for the number of generators
by degree is given by
\begin{equation}
\frac{1}{1-G_k(t)} = H_k(t) 
\end{equation}
With $k=3$ and $4$, one finds
\begin{equation}
\begin{split}
G_3(t) =&
t + t^2 + 3\,t^3 + 7\,t^4 + 17\,t^5 + 41\,t^6 + 99\,t^7 + 239\,t^8 +577\,t^9\\
&\,\,+ 1393\,t^{10} + 3363\,t^{11} + 8119\,t^{12} + 19601\,t^{13} + \dots
\end{split}
\end{equation}
which is Sequence~A001333 of~\cite{Slo}.

\begin{equation}
\begin{split}
G_4(t) =&
t + t^2 + 3\,t^3 + 13\,t^4 + 47\,t^5 + 173\,t^6 + 639\,t^7 + 2357\,t^8 +
8695\,t^9\\
&\,\, + 32077\,t^{10} + 118335\,t^{11} + 436549\,t^{12} + 1610471\,t^{13} +
\dots
\end{split}
\end{equation}
which is Sequence~A084519 of~\cite{Slo}.
With $k=5$, the sequence is not (yet) in~\cite{Slo}.

\begin{theorem}
\label{thm-subHopf}
$\DSym^{(k)}$ is a Hopf subalgebra of $\FQSym$.
\end{theorem}

\Proof
We already know that $\DSym^{(k)}$ is a subalgebra of $\FQSym$. So there
only remains to prove that $\DSym^{(k)}$ is a subcoalgebra of $\FQSym$.

Thanks to the definition of the coproduct of $\G_\sigma$ in $\FQSym$, this
amounts to a trivial combinatorial property:
if two permutations of size $n$ have same standardized word on each factor of
a given size, then, adding given letters all greater than $n$ (or all smaller
than $1$) at the same position in both permutations does not change the
property: the resulting permutations also have same standardized word on each
factor of the previous size.
\qed

\begin{corollary}
\label{thm-inclusions}
$\DSym^{(k)}$ is a Hopf subalgebra of $\DSym^{(l)}$ for $k<l$.
In particular, the $\DSym^{(k)}$ interpolate between $\NCSF=\DSym^{(2)}$ and
$\FQSym=\DSym^{\infty}$.
\qed
\end{corollary}

\subsection{Duality}

Let us denote by $\DQSym^{(k)}={\DSym^{(k)}}^*$, the dual bialgebra of
$\DSym^{(k)}$.
Dualizing Theorem~\ref{thm-subHopf}, we obtain:
\begin{theorem}
$\DQSym^{(k)}={\DSym^{(k)}}^* := \FQSym / \equiv_k$ is a noncommutative
(for $k>2$) and non-cocommutative Hopf algebra.
It is also a (non-free) dendriform quotient of $\FQSym$.
\end{theorem}

We shall write $\kF_C=\overline{\F_{\sigma}}$ where $C$ is the $k$-descent
code of $\sigma$.
For example, $52413\equiv_3 21543$, so that $42531$ and $21543$ have same
$3$-descent composition. We have
\begin{equation}
\begin{split}
\Delta \F_{42531} =& \
\F_{42531}\otimes 1 + 
\F_{3142}\otimes \F_{1} +
\F_{213}\otimes \F_{21}\\& +
\F_{21}\otimes \F_{321} +
\F_{1}\otimes \F_{2431} +
1\otimes \F_{42531},
\end{split}
\end{equation}
and
\begin{equation}
\begin{split}
\Delta \F_{21543} =& \
\F_{21543}\otimes 1 + 
\F_{2143}\otimes \F_{1} +
\F_{213}\otimes \F_{21}\\& +
\F_{21}\otimes \F_{321} +
\F_{1}\otimes \F_{1432} +
1\otimes \F_{21543}.
\end{split}
\end{equation}

Moreover,
\begin{equation}
\F_{42531} \F_{1}  =
\F_{425316} +
\F_{425361} +
\F_{425631} + 
\F_{426531} +
\F_{462531} +
\F_{642531},
\end{equation}
and
\begin{equation}
\F_{21543} \F_{1}  =
\F_{215436} +
\F_{215463} + 
\F_{215643} +
\F_{216543} + 
\F_{261543} + 
\F_{621543},
\end{equation}
and one easily checks that the indices of both expressions match in a
one-to-one correspondence by the $\equiv_3$ relation on the inverse
permutations.



\footnotesize
\begin{thebibliography}{aa}
\bibitem{NCSF7}{\sc G. Duchamp, F. Hivert, J.-C. Novelli},
and {\sc J.-Y. Thibon},
{\it Noncommutative Symmetric Functions VII: Free Quasi-Symmetric Functions
Revisited}, preprint, math.CO/0809.4479.
%
\bibitem{NCSF6}{\sc G. Duchamp, F. Hivert}, and {\sc J.-Y. Thibon},
{\it Noncommutative symmetric functions VI: free quasi-symmetric functions and
related algebras},
Internat. J. Alg. Comput. {\bf 12} (2002), 671--717.
%
\bibitem{NCSF1}{\sc I.M. Gelfand, D. Krob, A. Lascoux, B. Leclerc,
V.~S. Retakh}, and {\sc J.-Y. Thibon},
{\it Noncommutative symmetric functions},
Adv. in Math. {\bf 112} (1995), 218--348.
%
\bibitem{HTm}{\sc F. Hivert} and {\sc N. Thi\'ery},
{\it MuPAD-Combinat, an open-source package for research in algebraic
combinatorics},
S\'em. Lothar.  Combin. {\bf 51} (2004), 70p. (electronic).
%
\bibitem{NCSF2}{\sc D. Krob, B. Leclerc}, and {\sc J.-Y. Thibon},
{\it Noncommutative symmetric functions II: Transformations of alphabets},
Internal J. Alg. Comput. {\bf 7} (1997), 181--264.
%
\bibitem{NCSF4}{\sc D. Krob} and {\sc J.-Y. Thibon},
{\it Noncommutative symmetric functions IV{\,}: Quantum linear groups and
Hecke algebras at $q=0$}, J. Alg.  Comb.{\bf 6} (1997), 339--376.
%
\bibitem{Mcd}{\sc I.G. Macdonald},
{\it Symmetric functions and Hall polynomials},
2nd ed., Oxford University Press, 1995.
%
\bibitem{MR1}{\sc C.~Malvenuto} and {\sc C.~Reutenauer},
{\it Duality between quasi-symmetric functions and Solomon descent algebra},
J. Algebra, {\bf 177} (1995), 967--982.
%
\bibitem{Nov}{\sc J.-C. Novelli},
{\it On the hypoplactic monoid},
Disc. Math. {\bf 217} (2000), 315--336.
%
\bibitem{Reu}{\sc C. Reutenauer},
{\it Free Lie algebras},
Oxford University Press, 1993.
%
\bibitem{Slo}{\sc N.J.A. Sloane},
{\it The On-Line Encyclopedia of Integer Sequences},\\
\verb+http://www.research.att.com/~njas/sequences/+
%
\bibitem{So} \sc L. Solomon, \it A Mackey formula in the group
ring of a Coxeter group\rm, J. Algebra, {\bf 41}, (1976), 255-268.

%
\end{thebibliography}
\end{document}